\documentclass[leqno,12pt]{amsart}

\usepackage{palatino}
\usepackage[mathcal]{euler}
\usepackage[margin=1.5in]{geometry}
\usepackage{microtype}

\usepackage{amssymb,amsmath,amsthm,enumerate,slashed}
\usepackage[latin1]{inputenc}
\usepackage[T1]{fontenc}
\makeatletter
\input{ot1ppl.fd}
\DeclareFontShape{OT1}{ppl}{m}{up}{<->ssub*ppl/m/n}{}
\makeatother
\usepackage{mathrsfs}
\usepackage[dvipsnames,svgnames,table]{xcolor}
\usepackage{comment}

\definecolor{DarkBlue}{rgb}{0.2,0.2,0.4}
\usepackage[colorlinks=true, linkcolor=DarkBlue, urlcolor= DarkBlue, citecolor=DarkBlue]{hyperref}

\usepackage{graphicx}
\usepackage{epstopdf,epsfig}

\usepackage{xy}
\xyoption{all}
\usepackage{marginnote}
\usepackage[mmddyyyy]{datetime}

\usepackage[normalem]{ulem}

\makeatletter
\def\l@subsection{\@tocline{2}{0pt}{2.5pc}{5pc}{}}
\makeatother

\newcommand{\B}{\mathcal{B}}
\newcommand{\Ccal}{\mathcal{C}}

\newcommand{\norm}[1]{\left\|#1\right\|}

\DeclareMathOperator{\Tr}{Tr}

\DeclareMathOperator{\ran}{ran}
\DeclareMathOperator{\rank}{rank}

\DeclareMathOperator{\Span}{span}

\swapnumbers
\numberwithin{equation}{section}

\theoremstyle{plain}
\newtheorem*{theorem*}{Theorem}
\newtheorem*{lemma*}{Lemma}

\newtheorem{theorem}{Theorem}[section]

\newtheorem{lemma}[theorem]{Lemma}
\newtheorem{proposition}[theorem]{Proposition}

\theoremstyle{definition}
\newtheorem{construction}[theorem]{Construction}
\newtheorem{definition}[theorem]{Definition}
\newtheorem*{definition*}{Definition}

\newtheorem*{notation*}{Notation}

\theoremstyle{remark}
\newtheorem{remark}[theorem]{Remark}

\newtheorem*{remark*}{Remark}

\def\clap#1{\hbox to 0pt{\hss#1\hss}}

\usepackage{cancel} 

\title[Quasicentral modulus for almost normal operators]{Quasicentral modulus for almost normal operators}
\author{Shiqi Liu}
\address{S. Liu, University of New South Wales, Kensington, NSW, Australia}
\email{liushiqi127@gmail.com}
\author{Hongyin Zhao}
\address{H. Zhao, Central South University, Changsha, Hunan, China}
\email{hongyin.zhao@csu.edu.cn}
\thanks{Shiqi Liu is supported by the ARC grant DP230100434.
Hongyin Zhao is supported by China Postdoctoral Science Foundation -- Hunan Joint Support Program, grant 2025T001HN}

\begin{document}

\begin{abstract}
	In 1984, Voiculescu conjectured that the Hilbert--Schmidt quasicentral modulus vanishes for every almost normal operator. We prove this conjecture.
\end{abstract}

\maketitle

\section{Introduction}

Let $H$ be a separable complex Hilbert space. A bounded operator
$T{\in}\B(H)$ is \emph{almost normal} if
$[T^*,T]=T^*T-T\, T^*$ belongs to the trace class $\Ccal_1(H)$.
Write $\Ccal_2(H)$ for the Hilbert--Schmidt class.

Voiculescu introduced the quasicentral modulus for norm ideals
in \cite{Voi79}. We use the following sequential definition
\cite[Props.~1.1--1.2]{Voi79}.

\begin{definition}[Quasicentral Modulus for the Hilbert--Schmidt class]
	\label{def:modulus}
For $T\in \B(H)$, the \emph{quasicentral modulus} for the Hilbert--Schmidt class is
\begin{equation}\label{eq:modulus}
	k_2(T)=\inf_{(R_n)}\liminf_{n\to\infty}\norm{[T,R_n]}_2,
\end{equation}
where the infimum is taken over all sequences $(R_n)$ of finite-rank
operators satisfying $0\leqslant R_n\leqslant I$ and $R_n\to I$
strongly.
\end{definition}

The purpose of this paper is to prove the following theorem.
\begin{theorem}\label{thm:main}
If $T\in\B(H)$ and $[T^*,T]\in\Ccal_1(H)$, then $k_2(T)=0$.
\end{theorem}

Voiculescu conjectured this result in 1984 \cite{Voi84}.
He restated it in 2014 \cite[Sec.2.7]{Voi14}
and 2016 \cite[Sec.5.4]{Voi16}.
We prove it by constructing a sequence of finite-rank positive
contractions $F_{X_n}$ such that
\begin{equation}\label{eq:approximate-identity}
 F_{X_n}\longrightarrow I\quad\text{strongly},\qquad
 \norm{[T,F_{X_n}]}_2\longrightarrow0.
\end{equation}

To do so, we use Brown's finite-rank trace theorem and
Stone's formula for graph projections.
Our constructions and graph estimates hold
for every $T\in\B(H)$.
The trace-class condition is used only in the final
dominated convergence argument.

\section{Finite-rank interpolation}

In 1994, Brown proved a trace-vanishing theorem for finite-rank
commutators under Schatten-class assumptions
\cite[Lemma and Theorem, pp.~171--172]{Brown94}.
We need only the case $p=\infty$, $q=2$ in Brown's notation:
one operator is bounded and the other is Hilbert--Schmidt.
For completeness, we include a proof of this special case.

\begin{lemma}[Brown]\label{lem:trace}
Let $T\in\B(H)$ and $Y\in\Ccal_2(H)$. If $[T,Y]$ has finite rank,
then $\Tr[T,Y]=0$.
\end{lemma}

\begin{proof}
Unless otherwise stated, all identities and estimates below
hold for every integer $n\geqslant0$.
Put $F=[T,Y]$. We may suppose $F\ne0$.
Let $d=\rank F$ and, for $n\geqslant0$, set
\[
 M_n=\Span\{T^\alpha Y^\beta h:h\in\ran F,\ \alpha,\beta\in\mathbb Z_{\geqslant0},\ \alpha+\beta\leqslant n\}.
\]
For $n\geqslant 0$, let $P_n$ be the orthogonal projection
onto $M_n$ and set $E_n=P_n-P_{n-1}$, where $P_{-1}=0$.
In particular, $E_0=P_0$ is the orthogonal projection
onto $M_0=\ran F$.

The projections $E_n$ are mutually orthogonal.
For $n=0$, $\rank E_0=d$.
For $n{\geqslant}1$, $M_n$ is spanned by $M_{n-1}$ and $n{+}1$
subspaces
\[
 T^\alpha Y^{n-\alpha}(\ran F),\qquad
 0\leqslant\alpha\leqslant n,
\]
each of dimension at most $d$.
Hence,
\begin{equation}\label{eq:layer-rank}
	\rank E_n\leqslant d(n+1).
\end{equation}

We show the inclusions
\begin{equation}\label{eq-inclusion-fr}
	TM_n \subseteq M_{n+1}, \quad YM_n \subseteq M_{n+1}.
\end{equation}
The first inclusion is by definition. For the second inclusion, we use the algebraic identity
\[
 YT^\alpha=T^\alpha Y-\sum_{\ell=0}^{\alpha-1}T^\ell F T^{\alpha-1-\ell},
 \qquad \alpha \in\mathbb Z_{\geqslant0}.
\]
Applied to $Y^\beta h$ with $\alpha{+}\beta {\leqslant} n$, 
the first term belongs to $M_{n+1}$;
since $\ell{\leqslant} \alpha {-} 1 {<} n$, 
every remaining term belongs to $T^\ell(\ran F)\subseteq M_n$.

For $Z=T$ or $Y$, the inclusions in
\eqref{eq-inclusion-fr} give
\[
 (I-P_n)ZP_{n-1}=0,\quad
 P_{n+1}ZE_n=ZE_n.
\]
When $n=0$, the first equality is given by $P_{-1}=0$.
Using $E_n=P_n-P_{n-1}$, we can write the off-diagonal
block of $Z$ as
\begin{equation}\label{eq:crossing}
 (I-P_n)ZP_n=(I-P_n)ZE_n=E_{n+1}ZE_n.
\end{equation}

Since $\ran F\subseteq M_n$, we have $P_nF=F$.
Trace cyclicity gives
\begin{equation}\label{eq-tr-compress}
 \Tr(P_nFP_n)=\Tr(P_n^2F)=\Tr F,
 \qquad n\geqslant0.
\end{equation}
Inserting $I=P_n+(I-P_n)$ between the factors in $F=TY-YT$ gives
\[
 P_nFP_n=[P_nTP_n,P_nYP_n]
     +P_nT(I-P_n)YP_n-P_nY(I-P_n)TP_n.
\]
Since $P_n$ has finite rank, the commutator term
has trace zero by cyclicity of the trace.
Taking traces and using \eqref{eq:crossing} and
\eqref{eq-tr-compress}, we obtain
\begin{align}\label{eq:boundary-trace}
 \Tr F
 &=\Tr(P_nTE_{n+1}YE_n)-\Tr(P_nYE_{n+1}TE_n)\nonumber \\
 &=\Tr(E_nTE_{n+1}YE_n)-\Tr(E_nYE_{n+1}TE_n).
\end{align}
The second equality is due to $E_n^2=E_n$, $E_nP_n=E_n$, and trace cyclicity.

Recall that the Hilbert--Schmidt ideal $\Ccal_2$ is a Hilbert space,
 with inner product $\langle X,Y\rangle_2=\Tr(X^*Y)$.
By Cauchy--Schwarz and the ideal property,
\eqref{eq:boundary-trace} gives
\begin{align*}
 |\Tr F|
 &\leqslant
 \bigl|\langle(E_nTE_{n+1})^*,E_{n+1}YE_n\rangle_2\bigr|
 + \bigl|\langle(E_nYE_{n+1})^*,E_{n+1}TE_n\rangle_2\bigr|\\
 &\leqslant \norm{T}\norm{E_n}_2
 \bigl(\norm{E_{n+1}YE_n}_2+\norm{E_nYE_{n+1}}_2\bigr)\\
 &\leqslant \sqrt{d(n+1)}\norm{T}
 \bigl(\norm{E_{n+1}YE_n}_2+\norm{E_nYE_{n+1}}_2\bigr).
\end{align*}
The last inequality follows from
$\norm{E_n}_2^2=\rank E_n\leqslant d(n+1)$.
We square the estimate, divide by $n+1$, and sum.
Using $(s+t)^2\leqslant2(s^2+t^2)$, we obtain, for every $N\geqslant0$,
\[
\begin{aligned}
 |\Tr F|^2\sum_{n=0}^{N}\frac1{n+1}
 &\leqslant2d\norm{T}^2\sum_{n=0}^{N}
 \left(\norm{E_{n+1}YE_n}_2^2+\norm{E_nYE_{n+1}}_2^2\right)\\
 &\leqslant2d\norm{T}^2\norm{Y}_2^2.
\end{aligned}
\]
The blocks of $Y$ appearing in the sum are mutually orthogonal
in the Hilbert--Schmidt inner product, which gives the last inequality.
The right-hand side is independent of $N$, while the sum on the left diverges when $N\to \infty$. Thus $\Tr F=0$.
\end{proof}

Next, we use Brown's trace identity to obtain the finite-rank
operators with $[T, X]$ being small in Hilbert--Schmidt norm.

\begin{lemma}\label{lem:interpolation}
Let $T\in\B(H)$, let $P$ be a finite-rank orthogonal projection,
and let $\varepsilon>0$. There is a finite-rank operator $X$ such that
\[
 \norm{[T,X]}_2<\varepsilon,\qquad
 \norm{(X-I)P}_2<\varepsilon,\qquad
 \norm{P(X-I)}_2<\varepsilon.
\]
\end{lemma}

\begin{proof}
Consider the bounded complex-linear map
\[
 L_P:\Ccal_2(H)\longrightarrow\Ccal_2(H)^3,\qquad
 L_P(X)=([T,X],XP,PX).
\]
We show that $b=(0,P,P)$ lies in the closure of its range.

If $(Y,U,V)\in (\ran L_P)^{\perp}$, then, for every
$X\in\Ccal_2(H)$, trace cyclicity gives
\[
 \begin{aligned}
 0&=\Tr(Y^*[T,X])+\Tr(U^*XP)+\Tr(V^*PX)\\
  &=\Tr\bigl(([T^*,Y]+UP+PV)^*X\bigr).
 \end{aligned}
\]
Consequently $[T^*,Y]=-(UP+PV)$ has finite rank.
Lemma~\ref{lem:trace} gives \[\Tr(UP+PV)=0.\] 
So
\[
 \langle(Y,U,V),b\rangle
   =\Tr(U^*P)+\Tr(V^*P)
   =\overline{\Tr(UP+PV)}=0.
\]
Thus, $(0,P,P)\in (\ran L_P)^{\perp \perp} = \overline{\ran L_P}$.
Finally, finite-rank operators are dense in $\Ccal_2(H)$.
Continuity of $L_P$ therefore lets us approximate $b$ by $L_P(X)$
with $X$ finite rank, proving the three inequalities.
\end{proof}

\section{Graph coordinates and Commutator Estimates}\label{sec-3}

In 1951, Stone studied operators through the orthogonal projections
onto their graphs \cite{Stone51}.
Using this perspective, we associate to each finite-rank operator
$X$ the operators
\begin{equation}\label{eq-def-FG}
 F_X=XX^*(I+XX^*)^{-1},\qquad
 G_X=X^*X(I+X^*X)^{-1}.
\end{equation}
Both are finite-rank positive contractions.
For any $T\in\B(H)$, we estimate $\norm{[T,F_X]}_2^2$
in terms of $\norm{[T,X]}_2^2$ with a trace correction.
This construction transforms the finite-rank operators supplied
by Lemma~\ref{lem:interpolation} into finite-rank positive
contractions, as required by the definition of $k_2(T)$
(Definition~\ref{def:modulus}).

\subsection{Graph Coordinates and Projection}

\begin{lemma}\label{lem:graph-perp}
Let $X\in\B(H)$ and write
$\operatorname{Graph}(X)=\{(h,Xh): h\in H\}$. Then
\[
 \operatorname{Graph}(X)^\perp
 =\left\{(-X^*v,v):v\in H \right\}.
\]
\end{lemma}

\begin{proof}
A vector $(u,v)\in H \oplus H$ is orthogonal to $\operatorname{Graph}(X)$
if and only if $u+X^*v=0$.
\end{proof}

The maps
$h\mapsto(h,Xh)$ and $k\mapsto(-X^*k,k)$ identify $H$ with
$\operatorname{Graph}(X)$ and $\operatorname{Graph}(X)^\perp$,
respectively.
To make these identifications isometric, set
\begin{equation}\label{eq-AB}
	 A=(I+XX^*)^{-1/2},\qquad B=(I+X^*X)^{-1/2},
\end{equation}
and define $V,W:H\to H\oplus H$ by
\begin{equation}\label{eq:graph-isometries}
 Vh=(Bh,XBh),\qquad Wk=(-X^*Ak,Ak).
\end{equation}
A direct calculation shows that $V$ and $W$ preserve norms,
so they are isometries. Hence $V^*V=W^*W=I$, and
$V\,V^*$ and $WW^*$ are the orthogonal projections onto their
ranges. Since $A$ and $B$ are invertible, $V$ and $W$ are
isometric isomorphisms from $H$ onto $\operatorname{Graph}(X)$
and $\operatorname{Graph}(X)^{\perp}$, respectively.
The orthogonality of these ranges gives $V^*W=W^*V=0$.

\begin{definition}[Graph coordinates]\label{def:graph-coordinates}
Together, $V$ and $W$ give the following change of coordinates:
\[
\begin{array}{r@{\quad}c@{\;}c@{\;}c@{\;}c@{\;}c}
 U:
 & H\oplus H
 & \xrightarrow[\cong]{V\oplus W}
 & \operatorname{Graph}(X)\oplus\operatorname{Graph}(X)^\perp
 & \longrightarrow
 & H\oplus H\\[3pt]
 & (h,k)
 & \longmapsto
 & (Vh,Wk)
 & \longmapsto
 & Vh+Wk.
\end{array}
\]
We call $(h,k)$ the \emph{graph coordinates} of the vector
$Vh+Wk$. The map $U$ sends graph coordinates to the usual
direct-sum coordinates in $H\oplus H$.
\end{definition}

In matrix notation, $U=(V\;W)$, with block columns $V$ and $W$.
Using the identities above, we obtain
\[
\begin{aligned}
 U^*U
 &=\begin{pmatrix}V^*\\W^*\end{pmatrix}(V\;W)
 =I,\\
 UU^*
 &=(V\;W)\begin{pmatrix}V^*\\W^*\end{pmatrix}
 =V\,V^*+W\,W^*=I.
\end{aligned}
\]
Thus $U$ is unitary. The block rows of $U^*$ are
\begin{equation}\label{eq:graph-adjoints}
 V^*=\begin{pmatrix}B&BX^*\end{pmatrix},
 \qquad
 W^*=\begin{pmatrix}-AX&A\end{pmatrix}.
\end{equation}

We list two useful identities. The definitions in
\eqref{eq-def-FG} and \eqref{eq-AB}, and the relation
$(I+XX^*)X=X(I+X^*X)$, give
\begin{equation}\label{eq:graph-identities}
 XB^2X^*=F_X,\qquad I-B^2=G_X.
\end{equation}

Let $\Pi_X$ denote the orthogonal projection of $H\oplus H$
onto $\operatorname{Graph}(X)$.
Recall that $V$ is an isometry with range $\operatorname{Graph}(X)$,
so $\Pi_X=V\, V^*$.
In the usual coordinates,
\begin{equation}\label{eq:graph-corners}
	\Pi_X=V\, V^*
	= \begin{pmatrix} B \\ XB \end{pmatrix}
		 \begin{pmatrix} B & BX^* \end{pmatrix}
	=\begin{pmatrix}B^2&B^2X^*\\XB^2&XB^2X^*\end{pmatrix}
	=\begin{pmatrix}I-G_X&B^2X^*\\XB^2&F_X\end{pmatrix}.
\end{equation}
Stone called this block matrix the \emph{characteristic matrix}
of $X$ and gave its explicit formula \cite[Def.1, Thm.4]{Stone51}.

In graph coordinates, $\Pi_X$ acts by $(h,k)\mapsto(h,0)$.
Indeed,
\begin{equation}\label{eq:graph-diagonal}
 U^*\Pi_XU
 = \begin{pmatrix} V^* \\ W^* \end{pmatrix}
 	\cdot V\,V^* \cdot
 	\begin{pmatrix} V & W \end{pmatrix}
 =\begin{pmatrix}I&0\\0&0\end{pmatrix}.
\end{equation}


\subsection{Commutator Estimates}

In the usual coordinates, the characteristic matrix formula
\eqref{eq:graph-corners} shows that $[T,F_X]$ is the lower-right
block of $[T\oplus T,\Pi_X]$.
This gives the following Hilbert--Schmidt norm estimate.

\begin{lemma}\label{lem:graph-corner}
Let $T\in\B(H)$ and let $X$ be a finite-rank operator. 
Then $[T\oplus T,\Pi_X]$ has finite rank, and
\begin{equation}\label{eq:corner}
	\norm{[T,F_X]}_2\leqslant\norm{[T\oplus T,\Pi_X]}_2.
\end{equation}
\end{lemma}

\begin{proof}
By \eqref{eq:graph-corners},
\[
 [T\oplus T,\Pi_X]
 =\begin{pmatrix}
 -[T,G_X]&[T,B^2X^*]\\
 [T,XB^2]&[T,F_X]
 \end{pmatrix}.
\]
Since $X$ has finite rank, so do $F_X$, $G_X$, $B^2X^*$
and $XB^2$. Hence every block has finite rank.
The squared Hilbert--Schmidt norm of an operator is
the sum of the squared Hilbert--Schmidt norms of its blocks,
which proves \eqref{eq:corner}.
\end{proof}

We next compare $L\in \B(H)$ with the graph compression $V^*(L\oplus L)V$.

\begin{lemma}\label{lem:graph-trace}
Let $L\in\B(H)$ and let $X$ be a finite-rank operator.
Then $V^*(L\oplus L)V-L$ has finite rank, and
\[
 \Tr\bigl(V^*(L\oplus L)V-L\bigr)
 =\Tr\bigl(L(F_X-G_X)\bigr).
\]
\end{lemma}

\begin{samepage}
\begin{proof}
By \eqref{eq:graph-identities}, $I-B^2=G_X$.
Since $G_X$ has finite rank and $I+B$ is invertible,
$B-I=-G_X(I+B)^{-1}$ has finite rank. Thus
\begin{equation*}
	V^*(L\oplus L)V-L = (B-I)LB+L(B-I)+BX^*LXB
\end{equation*}
has finite rank. Cyclicity of the trace and
\eqref{eq:graph-identities} give
\begin{align*}
 \Tr\bigl(V^*(L\oplus L)V-L\bigr)
 & =\Tr\left(L(B^2-I+XB^2X^*)\right)\\
 & =\Tr\bigl(L(F_X-G_X)\bigr).\qedhere
\end{align*}
\end{proof}
\end{samepage}

\begin{samepage}
Using graph coordinates and the trace formula above, we obtain further results as follows.
\begin{proposition}\label{lem:graph}
Let $T\in\B(H)$ and let $X$ be a finite-rank operator. 
Then
\begin{multline}\label{eq:graph-identity}
 \norm{[T\oplus T,\Pi_X]}_2^2\\
 	=2\norm{(I+XX^*)^{-1/2}[T,X](I+X^*X)^{-1/2}}_2^2
 		+\Tr\bigl([T^*,T](G_X-F_X)\bigr).
\end{multline}
In particular,
\begin{equation}\label{eq:graph-estimate}
\norm{[T,F_X]}_2^2
 \leqslant 2\norm{[T,X]}_2^2
 +\bigl|\Tr\bigl([T^*,T](G_X-F_X)\bigr)\bigr|.
\end{equation}
\end{proposition}
\end{samepage}

\begin{proof}
Put $S=T\oplus T$. In graph coordinates, write
\[
 U^*SU=\begin{pmatrix}a&b\\c&d\end{pmatrix},
 \qquad a=V^*SV,\quad c=W^*SV.
\]
Combining this with \eqref{eq:graph-diagonal}, we calculate
\begin{equation}\label{eq:offdiagonal}
 U^*[S,\Pi_X]U
 = \left[\begin{pmatrix} a & b \\ c & d \end{pmatrix},  \begin{pmatrix} I & 0 \\ 0 & 0 \end{pmatrix} \right]
 =\begin{pmatrix}0&-b\\c&0\end{pmatrix}.
\end{equation}
By Lemma~\ref{lem:graph-corner}, $[S,\Pi_X]$ has finite rank,
so \eqref{eq:offdiagonal} shows that $b$ and $c$ have finite rank.
Since $U$ is unitary,
\begin{equation}\label{eq-bcc-norm}
	\norm{[S,\Pi_X]}_2^2=\norm{b}_2^2+\norm{c}_2^2.
\end{equation}

The formulas for $V$ in \eqref{eq:graph-isometries}
and $W^*$ in \eqref{eq:graph-adjoints} give
\begin{equation}\label{eq-graph-c}
 c=W^*SV
 	=\begin{pmatrix} -AX & A \end{pmatrix} 
 		\begin{pmatrix} T &  \\ & T \end{pmatrix}
 		\begin{pmatrix}  B \\ XB	\end{pmatrix}
 	=A[T,X]B.
\end{equation}
Since $A$ and $B$ are contractions by \eqref{eq-AB},
\begin{equation}\label{eq-c-norm-est}
	\norm{c}_2\leqslant\norm{[T,X]}_2.
\end{equation}

Now, we calculate $\norm{c}_2^2-\norm{b}_2^2$. Note that
\begin{equation*}
	a - T =V^*SV -T = V^*(T\oplus T)V -T.
\end{equation*}
Lemma~\ref{lem:graph-trace}, applied with $L=T$,
shows that $a-T$ has finite rank. Write
\[
 [a^*,a]-[T^*,T]
 =[a^*,a-T]+[a^*-T^*,T].
\]
Each commutator on the right has a finite-rank factor.
So the left-hand side is also finite rank. 
Cyclicity of the trace gives
\begin{equation}\label{eq:compression-trace}
	\Tr\bigl([a^*,a]-[T^*,T]\bigr)=0.
\end{equation}

We calculate $[S^*,S]$ in graph coordinates and take its
upper-left block. Subtracting $[T^*,T]$ gives
\begin{equation*}
 V^*[S^*,S]V-[T^*,T]
 =\bigl([a^*,a]-[T^*,T]\bigr)+c^*c-bb^*.
\end{equation*}

All terms on the right have finite rank, and the first has
trace zero by \eqref{eq:compression-trace}.
Since $[S^*,S]=[T^*,T]\oplus[T^*,T]$, taking traces and applying
Lemma~\ref{lem:graph-trace} with $L=[T^*,T]$ gives
 \begin{align}\label{eq-cb-difference}
 \norm{c}_2^2-\norm{b}_2^2
 &=\Tr\bigl(V^*[S^*,S]V-[T^*,T]\bigr) \nonumber \\
 &=\Tr\bigl([T^*,T](F_X-G_X)\bigr).
 \end{align}
Combining this identity with \eqref{eq-bcc-norm} gives
\[
 \norm{[T\oplus T,\Pi_X]}_2^2
 =2\norm{c}_2^2
 +\Tr\bigl([T^*,T](G_X-F_X)\bigr).
\]
Substituting \eqref{eq-graph-c} proves \eqref{eq:graph-identity}.
The estimate \eqref{eq:graph-estimate}
then follows from \eqref{eq:corner} and the bound
\eqref{eq-c-norm-est} on $\norm{c}_2$.
\end{proof}

\begin{remark}
	The estimate \eqref{eq-c-norm-est} controls $c$ in terms of $[T,X]$. 
	A direct estimate of $b=-B[T,X^*]A$ would involve $[T,X^*]$. 
	But $\norm{[T,X^*]}_2$ need not be small when $\norm{[T,X]}_2$ is small.
	However,  using Lemma~\ref{lem:graph-trace}, 
	we express $\norm{c}_2^2-\norm{b}_2^2$ as a trace involving $[T^*,T]$, see \eqref{eq-cb-difference}. 
\end{remark}

\section{Positive approximation}

\begin{construction}\label{con:scaled-approximants}
Let $T\in\B(H)$, and let $\{P_n\}_{n=1}^{\infty}$ be any increasing sequence of
finite-rank orthogonal projections converging strongly to $I$.
Put
\[
 D=\bigcup_n P_nH.
\]
Then $D$ is a dense subspace of $H$.
For any $n\geqslant1$, Lemma~\ref{lem:interpolation} gives a finite-rank
operator $A_n$ such that
\[
 \norm{[T,A_n]}_2<\frac1{n^2},\qquad
 \norm{(A_n-I)P_n}_2<\frac1{n^2},\qquad
 \norm{P_n(A_n-I)}_2<\frac1{n^2}.
\]
For every $h\in D$, we have $P_nh=h$ for all sufficiently
large $n$. Hence these inequalities give $A_nh\to h$ and
$A_n^*h\to h$.
Set
\[
 X_n=nA_n,
\]
thus
\begin{equation}\label{eq:scaled-commutator}
 \norm{[T,X_n]}_2<\frac1n.
\end{equation}
For each $n$, let the graph projection $\Pi_{X_n}$ and the
finite-rank positive contractions $F_{X_n}$ and $G_{X_n}$
be as defined in Section~\ref{sec-3}.
\end{construction}

Berger used scaled graphs of exact Hilbert--Schmidt intertwiners
\cite[proof of Thm.~2, pp.~6--7]{Berger81}.
Here, we rescale the finite-rank operators $A_n$ supplied by
Lemma~\ref{lem:interpolation}.

\begin{lemma}\label{lem:graph-convergence}
For the sequence $X_n$ in Construction~\ref{con:scaled-approximants},
\[
 \Pi_{X_n}\longrightarrow\begin{pmatrix}0&0\\0&I\end{pmatrix}
 \quad\text{strongly on }H\oplus H.
\]
In particular,
\begin{equation}\label{eq:corner-limits}
 F_{X_n}\longrightarrow I,\qquad G_{X_n}\longrightarrow I
 \quad\text{strongly on }H.
\end{equation}
\end{lemma}

\begin{proof}
For $h\in D$, the definition of the graph and
Lemma~\ref{lem:graph-perp} give
\[
\begin{aligned}
 (h/n,A_nh)
 &\in\operatorname{Graph}(X_n)=\ran\Pi_{X_n},\\
 (A_n^*h,-h/n)
 &\in\operatorname{Graph}(X_n)^\perp=\ker\Pi_{X_n}.
\end{aligned}
\]
These vectors converge to $(0,h)$ and $(h,0)$, respectively.
The operator $I-\Pi_{X_n}$ vanishes on $\ran\Pi_{X_n}$,
and $\Pi_{X_n}$ vanishes on $\ker\Pi_{X_n}$.
Both are orthogonal projections and hence contractions, so
\begin{align*}
	\norm{(I-\Pi_{X_n})(0,h)}
	& =\norm{(I-\Pi_{X_n})\bigl((0,h)-(h/n,A_nh)\bigr)} \\
	& \qquad \qquad \qquad \leqslant\norm{(0,h)-(h/n,A_nh)}\longrightarrow0,\\
	\norm{\Pi_{X_n}(h,0)}
 	& =\norm{\Pi_{X_n}\bigl((h,0)-(A_n^*h,-h/n)\bigr)}\\
 	& \qquad \qquad \qquad \leqslant\norm{(h,0)-(A_n^*h,-h/n)}\longrightarrow0.
\end{align*}

Since $D\oplus D$ is dense in $H\oplus H$ and
$\norm{\Pi_{X_n}}\leqslant1$ for all $n{\geqslant} 1$, these limits give the claimed strong
convergence of $\Pi_{X_n}$. Taking the two diagonal blocks in
\eqref{eq:graph-corners} gives \eqref{eq:corner-limits}.
\end{proof}

\begin{proof}[Proof of Theorem~\ref{thm:main}]
For any $T\in\B(H)$,
Construction \ref{con:scaled-approximants} gives a sequence
of finite-rank positive contractions $F_{X_n}$.
By Lemma \ref{lem:graph-convergence}, $F_{X_n}{\to} I$ strongly.
We now show that \mbox{$\norm{[T,F_{X_n}]}_2 \to 0$}.

By \eqref{eq:graph-estimate} and \eqref{eq:scaled-commutator},
\begin{equation}\label{eq-main-pf-1}
 \norm{[T,F_{X_n}]}_2^2
 \leqslant\frac{2}{n^2}
 +\bigl|\Tr\bigl([T^*,T](G_{X_n}-F_{X_n})\bigr)\bigr|,
 \quad n\geqslant 1.
\end{equation}

It remains to show the trace term converges to $0$. 
Here, we use the assumption that $[T^*,T]$ is trace class.
Since it is compact and self-adjoint, it has an orthonormal
eigenbasis $\{e_j\}_{j=1}^{\infty}$ with corresponding eigenvalues $\lambda_j$.
Then 
\begin{equation}\label{eq:trace-limit}
  \Tr\bigl([T^*,T](G_{X_n}-F_{X_n})\bigr)
 	=\sum_j\lambda_j
   		\langle e_j,(G_{X_n}-F_{X_n})e_j\rangle.
\end{equation}

Since $F_{X_n}$ and $G_{X_n}$ are contractions,
$\norm{G_{X_n}-F_{X_n}}\leqslant 2$.
Thus, for every $n$, the summands in \eqref{eq:trace-limit} satisfy
\[
 |\lambda_j\langle e_j,(G_{X_n}-F_{X_n})e_j\rangle|
 \leqslant 2|\lambda_j|
 \quad\text{for all } j\geqslant 1.
\]
The trace-class assumption on $[T^*,T]$ gives
\[
 \sum_j 2|\lambda_j|=2\norm{[T^*,T]}_1<\infty.
\]
By Lemma~\ref{lem:graph-convergence},
$G_{X_n}-F_{X_n}\to0$ strongly. For each fixed $j$,
\[
 \lambda_j\langle e_j,(G_{X_n}-F_{X_n})e_j\rangle
 \to0,\qquad n\to+\infty.
\]
Applying dominated convergence to the series in
\eqref{eq:trace-limit}, we have
\[
 \Tr\bigl([T^*,T](G_{X_n}-F_{X_n})\bigr)
 \to0,\qquad n\to+\infty.
\]

Combining this limit with \eqref{eq-main-pf-1} proves
$\norm{[T,F_{X_n}]}_2\to0$.
Hence, $k_2(T)=0$ by Definition~\ref{def:modulus}.
\end{proof}

\begin{remark}
	Almost normality is used only in the final dominated convergence argument.
	The rest of the paper holds for every $T\in\B(H)$.
\end{remark} 


\section*{Acknowledgments}
We thank Fedor Sukochev for introducing us to Voiculescu's pioneering work and for his helpful comments on earlier drafts of this paper. 

\section*{Use of artificial intelligence}
AI assisted us in exploring and testing ideas. It also helped us review mathematical arguments, search the literature, and edit the manuscript.
The AI models we used were GPT Sol 6.1 and Astra 6. For AI agents, we used our own customized version of Danus and the original version of Pharos. 
We would like to thank the AI4Math team at Peking University for developing these two open-source agents.


\bibliography{ref.bib}
\bibliographystyle{alpha}
\end{document}